\documentclass[a4paper,11pt]{amsart}

\usepackage[english]{babel}
\usepackage{amsmath, enumerate, amsfonts, amssymb, amsthm}

\theoremstyle{plain}
\newtheorem{theo}{Theorem}[section]
\newtheorem{theoI}{Theorem}
\newtheorem{prop}[theo]{Proposition}
\newtheorem{lem}[theo]{Lemma}

\newtheorem{corI}[theoI]{Corollary}
\newtheorem{clai}[theo]{Claim}

\theoremstyle{definition}
\newtheorem{defin}[theo]{Def\mbox{}inition}

\newtheorem*{rem*}{Remark}

\newtheorem{ex}[theo]{Example}

\newcommand{\dps}{\displaystyle}

\newcommand{\C}{\mathbb{C}}
\newcommand{\N}{\mathbb{N}}
\newcommand{\Z}{\mathbb{Z}}

\renewcommand{\k}{\mathbf{k}}

\newcommand{\CC}{\mathcal{C}}
\newcommand{\FF}{\mathcal{F}}
\newcommand{\Frac}{\mathrm{Frac}}

\newcommand{\QQ}{\mathcal{Q}}
\newcommand{\modQ}{{\mathrm{mod}\QQ}}
\newcommand{\PP}{\mathcal{P}}

\newcommand{\GG}{\mathcal{G}}

\newcommand{\spec}{\mathrm{Spec}}

\renewcommand{\O}{\mathcal{O}} 

\newcommand{\D}{\mathcal{D}} 
\newcommand{\Dn}{\mathcal{D}_n}
\newcommand{\Dnz}{\mathcal{D}_n \langle z \rangle}

\newcommand{\FDn}{\hat{\mathcal{D}}_n}
\newcommand{\FDnk}{\hat{\mathcal{D}}_n(\k)}
\newcommand{\FDnkz}{\hat{\mathcal{D}}_n(\k)\langle z \rangle}

\newcommand{\FDnCCz}{\hat{\mathcal{D}}_n(\CC) \langle z \rangle}

\newcommand{\z}{\langle z \rangle}

\newcommand{\dxi}{\partial _{x_i}}

\newcommand{\dx}[1]{\partial _{x_{#1}}}

\newcommand{\ddx}{\partial _x}

\newcommand{\dxsur}[2]{\frac{\partial {#1}}{\partial x_{#2}}}

\newcommand{\gr}{\mathrm{gr}} 
\newcommand{\lm}{\mathrm{lm}} 
\newcommand{\lc}{\mathrm{lc}} 
\newcommand{\Supp}{\mathrm{Supp}} 
\newcommand{\New}{\mathrm{New}}
\newcommand{\Exp}{\mathrm{Exp}} 

\newcommand{\W}{\mathcal{W}}
\newcommand{\E}{\mathcal{E}}

\title{Gr\"obner fan for analytic $D$-modules with parameters}

\author{Rouchdi Bahloul}
\address{Laboratoire de Math\'ematiques,
Universit\'e de Versailles St-Quentin-en-Yvelines,
45 avenue des Etats-Unis - B\^atiment Fermat,
78035 Versailles,
France}
\email{bahloul@math.uvsq.fr}

\begin{document}


\begin{abstract}
This is the first part of a work dedicated to the study of
Bernstein-Sato polynomials for several analytic functions depending on
parameters. The main result of this part is a constructibility result
for the analytic Gr\"obner fan of a parametric ideal in the ring of
analytic differential operators. In this part, the main tool is the
notion of generic reduced standard basis.
\end{abstract}


\maketitle


\section*{Introduction}\label{sec:intro}
This is the first part of a work dedicated to the study of
Bernstein-Sato polynomials of several variables for analytic functions
depending on parameters (part 2 is \cite{part2}).
In the present paper, we focus on a major step of this study
and which constitutes a result interesting by its own. It is a
constructibility result for the analytic Gr\"obner fan of a parametric
ideal in the ring of analytic differential operators.

The Gr\"obner fan for polynomials was introduced by Mora
and Robbiano in 1988
\cite{mora-rob} (but earlier works by Lejeune-Jalabert and Teissier
\cite{lejeune-teissier} already contain analogous
constructions). Since then the Gr\"obner fan has found numerous
applications: e.g. the Gr\"obner walk in commutative algebra.
The ``$D$-modules version'' has been treated by Assi et
al. \cite{acg1, acg2} and Saito et al. \cite{sst}. The algebraic
version \cite{acg1, sst} has a nice application to
GKZ-hypergeometric differential systems (see \cite[Chapters 2,
3]{sst}). The analytic version \cite{acg2}
(see also \cite{btaka} for a significant extension)
made possible a complete
proof of the existence of Bernstein-Sato polynomials for
\emph{several analytic} functions (Bahloul \cite{compos}).

Based on \cite{compos} our goal is to give new constructive results
concerning Bernstein-Sato polynomials for several
analytic functions depending on parameters (see the second part of
this work \cite{part2}). A major step towards this goal consists in
studying
the Gr\"obner fan of an ideal depending on parameters in the ring
of analytic differential operators. As the main tool for this study,
we shall use parametric standard bases for analytic differential
ideals (in fact we shall work in a formal setting).
Parametric Gr\"obner bases (for polynomials) have been
extensively studied. A local version also exists: see e.g. Greuel and
Pfister \cite{greuel-p} and  Aschenbrenner \cite{aschen}.
An algebraic differential version  has been initiated by Oaku
\cite{oaku} whose work inspired Leykin \cite{leykin} and Walther
\cite{walther}. The analytic (or formal) differential version has been
introduced in Bahloul \cite{jmsj} without using reduced standard
bases. In the present paper we add as a supplementary tool the use of
parametric standard bases which are reduced in a sense we shall
define later. Indeed, since the Gr\"obner fan is described using
reduced standard bases, this is necessary.

Let us summarize: In the first section we introduce some notations
and recall some facts about the Gr\"obner fan, the division
theorem and (reduced) standard bases. Section 2 contains the
main tool: reduced generic standard bases. Section 3 is devoted to the
proof of the main result which we state now:

\begin{theoI}\label{theo:main}
Let $x=(x_1,\ldots,x_n)$ and $y=(y_1,\ldots,y_m)$ be two sets of
variables.
Let $I$ be an ideal in $\C\{x,y\}[\dx{1}, \ldots, \dx{n}]$
the ring of analytic differential operators with parameters $y$ (it is
also described as the germ of
ring of relative differential operators $\D_{\C^{n+m} / \C^m}$).
Let $\QQ$ be a prime ideal in $\C\{y\}$. There exists $h(y) \in
\C\{y\} \smallsetminus \QQ$ such that for any $y_0 \in V(\QQ)$
with $h(y_0) \ne 0$, the analytic Gr\"obner fan of
$I_{|y=y_0}$ is constant. Here $V(\QQ)$ means the zero
set of a representative of $\QQ$ on a small polydisc.
\end{theoI}

\begin{rem*}
We have a similar result concerning the global Gr\"obner fan for an
ideal in $D[y]$, where $D$ is the Weyl algebra over a field of
characteristic $0$. The proof is easier because every
process is finite. One can find a proof in \cite[Chap. 6]{these}.
\end{rem*}

\begin{corI}\label{corI}
There exists a finite stratification of $(\C^m, 0)=\cup W$
made of locally closed subsets such that the analytic Gr\"obner
fan of $I_{|y}$ is constant along each member $W$ of the partition.
\end{corI}

Let $\E$ be the common refinement of all the Gr\"obner
fans of $I_{|y=y_0}$, $y_0 \in (\C^m, 0)$ (by Cor. \ref{corI},
there is only a finite number of such fans).
Then $\E$ is the smallest fan with the following property:
For any $w$ in an open cone of $\E$ and $y_0 \in
(\C^m,0)$, the
graded ideal $\gr^{w}(I_{|y=y_0})$ is constant.
We call $\E$ the \emph{comprehensive Gr\"obner fan} of $I$
(following the terminology relative to Gr\"obner bases).

The idea of introducing this fan is due to N. Takayama whom I thank.\\
{\bf Acknowledgement.} The author was supported by the JSPS as a
postoctoral fellow at Kobe university when the results of
this work were elaborated.

\section{Recalls: Gr\"obner fan, divisions}

In this section, we recall the homogenized ring of (formal)
differential operators $\FDnkz$ ($\k$ is a field with characteristic
$0$). We then recall the definition of the Gr\"obner fan as in
\cite{acg2, btaka}. Finally we recall the division theorem in $\FDnkz$
together with the notion of (reduced) standard basis (see loc. cit.).

Let $x=(x_1,\ldots,x_n)$ be a system of variables, we write
$\ddx=(\dx{1}, \ldots, \dx{n})$ for the corresponding partial
differentials. Let $z$ be another variable. The ring $\FDnk$ is the
ring of differential operators with coefficients in $\k[[x]]$.
As we will work with arbitrary orders, we will need to work in a
homogenized (or graded) version of $\FDnk$.
The ring $\FDnkz$ is defined as the $\k[[x]]$-algebra generated by
$\ddx$ and $z$ where the only non trivial relations are:

$[\dxi, a]= \dxsur{a}{i} z$ for $i=1,\ldots,n$ and $a\in \k[[x]]$.

When we replace $\k[[x]]$ with $\C\{x\}$ we obtain and write $\Dnz$.
We see that the previous relation preserves the total degree in the
$\dxi$'s and $z$, which makes $\FDnkz$ a graded algebra for this
degree. Note that it is isomorphic to the Rees algebra of $\FDn$
associated with the filtration by the degree.

Given $P\in \FDnk$, $P=\sum c_{\alpha \beta} x^\alpha \ddx^\beta$
($\alpha, \beta \in \N^n$, $c_{\alpha \beta}\in \k$), we
define $h(P) \in \FDnkz$ as $h(P)=\sum c_{\alpha \beta} x^\alpha
\ddx^\beta z^{d-|\beta|}$ where $d=\deg(P)$ is the degree of $P$ in the
$\dxi$'s. Given a (left) ideal $I \subset \FDnk$, we define $h(I)$ as
the ideal of $\FDnkz$ generated by $\{h(P) | P\in I\}$.

Let $w=(u,v) \in \Z^{n+n}$. We consider it as a weight vector on the
variables $(x,\ddx)$. We define $\W=\{w| \forall i, u_i\le 0, u_i+v_i
\ge 0\}$ as the set of \emph{admissible weight vectors}. With $w\in
\W$, we can associate a natural filtration on $\FDnkz$ and a
graded ring $\gr^w(\FDnkz)$. Let $I$ be a given ideal in $\FDnk$ (an
analoguous construction holds for $\Dn$).
For $w,w'\in \W$, we write $w \sim w'$ when $\gr^w(\FDnkz)=
\gr^{w'}(\FDnkz)$ and $\gr^w(h(I))=\gr^{w'}(h(I))$. The partition
of $\W$ given by this relation is called the (formal) {\em open
  Gr\"obner fan of $I$} and denoted by $\E(h(I))$. It is a finite
collection of open convex polyhedral cones (see Assi et
al. \cite{acg2}). Denote by $\bar{\E}(h(I))$ the set of closures
of these open cones and call it the closed Gr\"obner fan of $I$.
In Bahloul, Takayama \cite{btaka}, we proved 
$\bar{\E}(h(I))$ is a polyhedral fan improving results of
Assi et al. \cite{acg2}.

It is easy to prove that the analytic Gr\"obner fan of $I \subset \Dn$
coincides with the (formal) Gr\"obner fan of $\FDn(\C) I$ (see
\cite{acg2}, \cite{btaka}). That is why we shall mainly work in a
formal setting.

Let us now deal with divisions and standard bases.
A total monomial order $\prec$ on $\N^{2n}$ (or equivalently on
the monomials $x^\alpha \xi^\beta$, $\xi_i$ being a comutative
variable corresponding to $\dxi$) is said to be \emph{admissible}
if $x_i \prec 1$ and $x_i \xi_i \succ 1$. For a given weight vector
$w$, we can define the order $\prec_w$ by refining $w$ by
$\prec$ (i.e. we first use $w$ and then tie-break with $\prec$).
Obviously if $w$ and $\prec$ are admissible then so is
$\prec_w$. Given $\prec$ admissible, we may define $\prec^h$ on
$\N^{2n+1}$: $(\alpha, \beta, k) \prec^h (\alpha', \beta', k')$ iff
$|\beta|+k < |\beta'| +k'$ or ($|\beta|+k = |\beta'| +k'$ and
$(\alpha, \beta) \prec (\alpha', \beta')$).

Let $P=\sum c_{\alpha \beta k} x^\alpha \ddx^\beta z^k$ be in
$\FDnkz$, we define the support $\Supp(P) \subset \N^{2n+1}$ as the
set of $(\alpha, \beta, k)$ with $c_{\alpha \beta k} \ne 0$. When $P
\ne 0$ we define its leading exponent $\exp_{\prec^h}(P)=
\max_{\prec^h}\Supp(P)$, leading coefficient $\lc_{\prec^h}(P) =
c_{\exp_{\prec^h}(P)}$ and leading monomial
$\lm_{\prec^h(P)}(P)= \lc_{\prec^h}(P) \cdot (x,\ddx,
z)^{\exp_{\prec^h}(P)}$. From now on, we shall omit the subscript
$\prec^h$ when the context is clear. Let us now recall the
division theorem.

Let $P_1, \ldots, P_r \in \FDnkz$ and $\prec$ an admissible order. With
$\prec^h$, we define the partition
$\Delta_1 \cup \cdots \cup \Delta_r \cup \bar{\Delta}$
of $\N^{2n+1}$ as
$\Delta_1= \exp(P_1)+\N^{2n+1}$ and for $j\ge 2$,
$\Delta_j=(\exp(P_j)+\N^{2n+1}) \smallsetminus \cup_{k=1}^{j-1}
\Delta_k$.

\begin{theo}[{Division theorem, see \cite[Th. 7]{acg2},
    \cite[Th. 3.1.1]{btaka}}] \label{theo:div}
For $P \in \FDnkz$, there exists a unique $(Q_1, \ldots, Q_r, R) \in
(\FDnkz)^{r+1}$ such that $P=\sum_j Q_j P_j +R$ and
\begin{itemize}
\item
for $j\ge 1$, $Q_j=0$ or $\Supp(Q_j) + \exp(P_j) \subset \Delta_j$,
\item
$R=0$ or $\Supp(R) \subset \bar{\Delta}$.
\end{itemize}
$R$ is called the remainder of the division. 
\end{theo}

As an easy consequence of the formal division process, we have:

\begin{lem}\label{lem:denom}
Let $\CC$ be a commutative integral ring and $\FF$ a field containing
$\CC$. Let $P, P_1,\ldots,P_r$ be in $\FDn(\CC)\z$. Let us
consider the division of $P$ by the $P_j$'s in $\FDn(\FF)\z$
w.r.t. $\prec^h$: $P=\sum Q_j P_j +R$. Then the coefficients of $R$ and
of the $Q_j$'s have the following form:
\[\frac{c}{\prod_{j=1}^r \lc(P_j)^{d_j}} \text{ where } c \in \CC, \,
d_j \in \N.\]
\end{lem}

For a (non zero) left ideal $J \subset \FDnkz$, we define the set of
exponents $\Exp(J)=\{\exp(P) | P\in J, P \ne 0\}$. A standard basis $G
\subset J$ of $J$ (with respect to $\prec^h$) is defined by the
relation: $\Exp(J)=\cup_{g \in G} (\exp(g)+\N^{2n+1})$. Thanks to the
division theorem, we have: $G$ is standard basis iff (for any $P$,
$P\in J \iff R=0$, where $R$ is the remainder of the division of $P$
by $G$). By noetherianity of $\N^{2n+1}$ a standard basis always
exists (Dickson lemma).

Let us end these preliminaries with the notion of a reduced standard
basis.

\begin{defin}
A $\prec^h$ standard basis $G=\{g_1,\ldots,g_r\}$ of $J\subset \FDnkz$
is said to be
\begin{itemize}
\item
minimal if for any $F\subset \N^{2n+1}$, we have:

$\Exp_{<_L^h}(J)=\bigcup_{e\in F} (e+\N^{2n+1}) \Rightarrow
\{\exp(g_1),\ldots,\exp(g_r)\} \subseteq F$.
\item
reduced if it is minimal and if for any $j=1,\ldots,r$, $\lc(g_j)=1$
and

$(\Supp(g_j)\smallsetminus \exp(g_j)) \subseteq (\N^{2n+1} \smallsetminus
\Exp(J))$.
\end{itemize}
\end{defin}

Given an ideal $J \subset \FDnkz$ and $\prec$ admissible, a
$\prec^h$-reduced standard basis exists and is unique.\\
Let us sketch the existence:\\
Let $G_0$ be any standard basis. By removing unnecessary elements we
may assume $G_0$ to be minimal. Set $G_0=\{g_j; 1\le j \le r\}$. For
any $j$, divide $g_j-\lm(g_j)$ by $G_0$ and denote by $r_j$ the
remainder. The set $\{(\lm(g_j)+r_j)/\lc(g_j); \, 1\le j\le r\}$ is
then the reduced standard basis of $J$.

\section{(Reduced) generic standard bases}

In \cite{jmsj}, we gave a definition of generic standard bases. A more
general definition is necessary if we want to deal with the property
of reduceness. However for the proofs, we shall cite \cite{jmsj},
because they are similar. The symbol $\prec$ shall denote an
admissible order.

Let $\CC$ be a commutative integral unitary ring 
(not necessarily noetherian) for which we denote
by $\FF$ the fraction field, by $\spec(\CC)$ the spectrum.
We assume also that for any $\PP \in \spec(\CC)$ and $n\in \Z$,
if $n\in \PP$ then $1 \in \PP$.
For any ideal $\mathcal{I}$ in $\CC$, we denote by
$V(\mathcal{I})=\{\PP \in \spec(\CC); \mathcal{I} \subset \PP\}$ the
zero set defined by $\mathcal{I}$.

For any $\PP$ in $\spec(\CC)$ and $c$ in $\CC$, denote by $[c]_\PP$
the class of $c$ in $\CC/\PP$ and by $(c)_\PP$ this class viewed in the
fraction field $\FF(\PP)=\Frac(\CC/\PP)$. (The condition above on $\CC$
ensures that $\FF(\PP)$ has characteristic $0$.) The element $(c)_\PP$ is
called the \emph{specialization of} $c$ \emph{to} $\PP$.

We naturally extend these notations to elements in $\FDnCCz$ and we
extend $(\cdot )_\PP$ to elements of $\FDn(\FF)\z$ for which the
denominators of the coefficients are in $\CC \smallsetminus \PP$, i.e.
$\FDn(\CC_\PP)\z$ where $\CC_\PP$ is the localization w.r.t. $\PP$.

Now, given an ideal $J \subset \FDnCCz$, we define the specialization
$(J)_\PP$ of $J$ to $\PP$ as the ideal of $\FDn(\FF(\PP))\z$
generated by all the $(P)_\PP$ with $P \in J$.

\subsection{Generic standard basis on an irreducible affine scheme}

\

Fix a prime ideal $\QQ$ in $\CC$.
Let us start with some notations. We denote by $\FDn(\QQ)\z$ the ideal
of $\FDnCCz$ made of elements with all their coefficients in $\QQ$.
For $h\in \CC$, we denote by $\CC[h^{-1}]$ the localization of $\CC$ w.r.t.
$h$. The ring $\FDn(\CC[h^{-1}])\z$ shall be seen as the subring of
$\FDn(\FF)\z$ made of elements with coefficients $\frac{c}{c'}$ such
that $c'$ is a power of $h$. In the latter, if all the $c$ are in
$\QQ$, we obtain an ideal denoted by $\FDn(\QQ[h^{-1}])\z$.
Finally, $\langle \QQ \rangle$ denotes the ideal of $\FDn(\FF)\z$ made of
elements with coefficients $\frac{c}{c'}$ such that $c\in \QQ$. Notice
that the latter is a priori different from $\FDn(\FF)\z \QQ$ since we
don't suppose $\CC$ to be noetherian.

Now for an element $P$ in $\FDnCCz \smallsetminus \FDn(\QQ)\z$ or more
generally in $\FDn(\CC_\QQ)\z \smallsetminus \langle \QQ \rangle$,
let us write $P=\sum \frac{c_{\alpha \beta k}}{c'_{\alpha \beta k}}
x^\alpha \ddx^\beta z^k$ with
$c'_{\alpha \beta k} \in \CC \smallsetminus \QQ$. Then denote by
$\exp^\modQ(P)$ the maximum (w.r.t. $\prec^h$) of the $(\alpha,\beta,
k)$ with $c_{\alpha \beta k} \notin \QQ$. This is the leading exponent
of $P$ modulo $\QQ$.
In the same way, we define the leading coefficient $\lc^\modQ(P)$ and
leading monomial $\lm^\modQ(P)$ modulo $\QQ$.

In general we have $\exp^\modQ(PQ)= \exp^\modQ(P)+\exp^\modQ(Q)$
as for the usual leading exponent. However, there are some differences
with the usual situation, for example the leading coefficient $\modQ$
of $PQ$ is not equal to the product of that of them. They are equal
only modulo $\QQ$ so we will have to be careful.

Now for an ideal $J\subset \FDnCCz \smallsetminus \FDn(\QQ)\z$, we
define: $\Exp^\modQ(J)=\{\exp^\modQ(P) | P\in J \smallsetminus
\FDn(\QQ)\z \}$.
This set is stable by sums in $\N^{2n+1}$ thus by Dickson lemma:
\begin{equation}\label{eq:gensb}
\begin{cases}
\exists \{g_1,\ldots,g_r\}\subset J \text{ such that }\\
\Exp^\modQ(J)=
\bigcup_j (\exp^\modQ(g_j)+ \N^n).
\end{cases}
\end{equation}
This shall be a generic standard basis of $J$ on $V(\QQ)$.
However, this is not the definition we will adopt. In fact in the next
paragraph we will define the notion of reduced generic standard basis
and it will not be in the ring $\FDnCCz$ so we need a more general
definition:

\begin{defin}\label{def:genSB}
A generic standard basis (gen.s.b for short) of $J$ on $V(\QQ)$ is a
couple $(\GG,h)$ where
\begin{itemize}
\item[(a)]
$h\in \CC \smallsetminus \QQ$,
\item[(b)]
$\GG$ is a finite set in the ideal $\FDn(\CC[h^{-1}])\z \cdot J$ and
for any $g\in \GG$ the numerator of $\lc^\modQ(g)$ divides $h$,
\item[(c)]
$\dps \Exp^\modQ(J)=\bigcup_{g\in \GG}(\exp^\modQ(g) +\N^{2n+1})$.
\end{itemize}
\end{defin}

Above in (\ref{eq:gensb}), $(\{g_1,\ldots,g_r\}, \prod_j
\lc^\modQ(g_j))$ is a gen.s.b of $J$ on $V(\QQ)$.
Thus, Def. \ref{def:genSB} makes sens.\\
Notice that another way to state (b) is: For any $\PP \in V(\QQ)
\smallsetminus V(h)$, the specialization $(g)_\PP$ is well defined and
belongs to $(J)_\PP$ and $\exp((g)_\PP)$ is equal to
$\exp^\modQ(g)$. Note 
that $V(\QQ) \smallsetminus V(h)$ is non empty since $h\notin \QQ$.


\begin{prop}[Division modulo $\QQ$]\label{prop:divmodQ}
Let $h\in \CC \smallsetminus \QQ$ and $g_1,\ldots,g_r$ be in
$\FDn(\CC[h^{-1}])\z \smallsetminus \FDn(\QQ[h^{-1}])\z$ such that
each $\lc(g_j)$ divides $h$.
Let $\Delta_1 \cup \cdots\cup \Delta_r \cup \bar{\Delta}$ be the
partition of $\N^{2n+1}$ associated with the $\exp^\modQ(g_j)$.
Then for any $P$ in $\FDn(\CC[h^{-1}])\z$, there exist
$q_1,\ldots,q_r, R, T \in \FDn(\FF)\z$ such that
\begin{itemize}
\item[(o)]
$P=\sum_j q_j g_j +R+T$,
\item[(i)]
$\Supp(q_j)+\exp^\modQ(g_j) \subset \Delta_j$ if $q_j\ne 0$,
\item[(ii)]
$\Supp(R) \subset \bar{\Delta}$ if $R \ne 0$,
\item[(iii)]
the $q_j$ and $R$ are in $\FDn(\CC[h^{-1}])\z$ and $T$ is in
$\FDn(\QQ[h^{-1}])\z$.
\end{itemize}
Moreover, $(q_1,\ldots,q_r,R)$ is unique modulo
$\FDn(\QQ[h^{-1}])\z$. 
We call $R$ the remainder $\modQ$ of this division.
\end{prop}

\begin{proof}[Sketch of proof]
Write $g_j=g_j^{(1)} - g_j^{(2)}$ with $g_j^{(2)} \in \langle \QQ
\rangle$ and $\exp (g_j^{(1)})=\exp^\modQ(g_j)$ then divide $P$ by the
$g_j^{(1)}$'s in $\FDn(\FF)\z$ as in theorem \ref{theo:div}:
$P=\sum_j q_j g_j^{(1)} +R$. 
We have, $P= \sum_j q_j g_j +R +T$ with $T=\sum_j q_j g_j^{(2)}$.
Conditions (i) and (ii) are satisfied by theorem \ref{theo:div}.
The third one is a direct consequence of lemma \ref{lem:denom}.
The last statement comes from the unicity in Th. \ref{theo:div}
after specializing to $\QQ$.
\end{proof}



The main result concerning generic standard basis is the following.

\begin{theo}\label{theo:genSB}
Let $(\GG,h)$ be a gen.s.b of $J \subset \FDn(\CC)\z$ on $V(\QQ)$.
Then for any $\PP \in V(\QQ) \smallsetminus V(h)$:
\begin{itemize}
\item[(i)] $(\GG)_\PP \subset (J)_\PP$,
\item[(ii)] $\dps \Exp((J)_\PP)=\bigcup_{g\in \GG} (\exp((g)_\PP)+\N^n)
=\Exp^\modQ(J)$.
\end{itemize}
\end{theo}
In other words, $(\GG)_\PP$ is a standard basis of $(J)_\PP$ for a
generic $\PP \in V(\QQ)$ and $\Exp((J)_\PP)$ is generically constant
and equal to $\Exp^\modQ(J)$.

\begin{proof}
Exactly the same as for \cite[Th. 3.7]{jmsj}
%
\end{proof}

\subsection{Reduced generic standard bases}

The next result shall concern the existence of \emph{the} reduced
generic standard basis on $V(\QQ)$ (in fact we shall see that it is
unique ``modulo $\QQ$'').

Let $J$ be an ideal in $\FDn(\CC)\z$ and $\QQ$ be a prime ideal in
$\CC$.

\begin{theo}[Def{}inition-Theorem]\label{theo:red.gen.sb}
\noindent
\begin{itemize}
\item
There exists a gen.s.b $(\GG, h)$ of $J$ on $V(\QQ)$ such that $(\GG)_\QQ$
is the reduced standard basis of $(J)_\QQ$. Such a $(\GG,h)$ is called a
reduced generic standard basis (red.gen.s.b) of $J$ on $V(\QQ)$.
\item
If $(\GG,h)$ is a red.gen.s.b on $V(\QQ)$ then 
for any $\PP \in V(\QQ) \smallsetminus V(h)$, $(\GG)_\PP$ is the reduced
standard basis of $(J)_\PP$.
\end{itemize}
\end{theo}

\noindent
Such a red.gen.s.b is unique ``modulo $\langle \QQ \rangle$''.
More precisely:
\begin{lem}
Let $(\GG,h)$ and $(\GG',h')$ be two red.gen.s.b of $J$ on $V(\QQ)$
then
\begin{itemize}
\item
their cardinality and the set of their leading exponents $\modQ$ are
equal,
\item
if $g\in \GG$ and $g'\in \GG'$ satisfy $\exp^\modQ(g)=\exp^\modQ(g')$
then $g-g'$ belongs to $\FDn(\QQ[(hh')^{-1}])\z$.
\end{itemize}
\end{lem}

\begin{proof}
The first statement is trivial by unicity of reduced standard bases.
For the second one, we have $g-g' \in \FDn(\CC[hh'^{-1}])\z$ and by
the same argument of unicity, $(g)_\QQ -(g')_\QQ=0$ thus $g-g'\in
\FDn(\QQ[hh'^{-1}])\z$.
\end{proof}

\begin{proof}[Proof of the Theorem]
For the first statement, let $(\GG_0, h)$ be any gen.s.b of $J$ on $V(\QQ)$.
Set $\GG_0=\{g_1,\ldots, g_r\}$. By removing the unnecessary elements, we
may assume that it is minimal. For any $j$ we may assume $\lc^\modQ(g_j)$
to be unitary. For any $j$, let $r_j$ be the remainder $\modQ$ of the
division modulo $\QQ$ of $g_j-\lm^\modQ(g_j)$ by $\GG_0$.
Set $\GG=\{\lm^\modQ(g_j)+r_j | j=1,\ldots,r\}$. It is easy to check that
$(\GG,h)$ is a red.gen.s.b.

Let us prove the second statement. Let $(\GG,h)$ be a red.gen.s.b.
First, we know that for any $\PP \in V(\QQ) \smallsetminus V(h)$,
$(\GG)_\PP$ is a standard basis of $(J)_\PP$.
Moreover it is minimal since $\Exp((J)_\PP)=\Exp((J)_\QQ)$ and
$\exp((g)_\PP)=\exp^\modQ(g)=\exp((g)_\QQ)$ for any $g\in \GG$.
The latter also implies that it is unitary. It just remains to prove
that it is reduced.
But this follows from the fact that $(\GG)_\QQ$ is reduced and that
for any $g\in \GG$, $\Supp((g)_\PP) \subset \Supp((g)_\QQ)$ (since $\QQ
\subset \PP$).
\end{proof}

The following example shows that reduced generic standard bases need
to be defined in some extension as in Def. \ref{def:genSB}.
\begin{ex}
Take $f(x_1,x_2, y)=y x_2 - x_1 x_2+ x_1$
in $\C[y][[x]]$ (it is a commutative situation for simplification).
Take an order such that the leading exponent
(in terms of $x$) is $(0,1)$ (i.e. corresponding to $x_2$). Take
$\QQ=(0)$ in $\C[y]$ then the red.gen.s.b of the ideal generated by
$f$ is $x_2+ x_1/y + x_1^2/y^2+ x_1^3/y^3+ \cdots$. 
\end{ex}

One can find more details and more results on (reduced) generic
standard bases in Bahloul \cite{gen_compSB}.

\section{Back to Gr\"obner fans}

In this section we shall prove Theorem \ref{theo:main}. First let us
recall two results needed for the proof.

Take $I \subset \FDnk$. For $w,w' \in \W$ we have defined an
equivalence relation $w \sim w'$ and defined the open Gr\"obner fan as
the collection of the equivalence classes. For $w\in \W$, let us
denote by $C_w(h(I))$ the equivalence class of $w$.

\begin{clai}[Recall 1] \label{recall1}
Let $w' \in C_w(h(I))$ and $\prec$ any admissible order
then the reduced standard bases of $h(I)
\subset \FDnkz$ with respect to $\prec_w^h$ and $\prec_{w'}^h$
coincide.
\end{clai}

This is by construction of the Gr\"obner fan because
$\Exp_{\prec_w^h}(h(I))=\Exp_{\prec_{w'}^h}(h(I))$ (see
\cite{acg2}). Another way to see this equality: By definition
$\gr^w(h(I))= \gr^{w'}(h(I))$, so $\Exp_{\prec_w^h}(h(I))=
\Exp_{\prec^h}(\gr^w(h(I))) = \Exp_{\prec^h}(\gr^{w'}(h(I))) =
\Exp_{\prec_{w'}^h}(h(I))$ (see \cite[Lemma 3.2.2]{btaka}).\\

For $g\in h(I)$, we define the Newton polyhedron as the following
convex hull: $\mathrm{New}(g)=
\mathrm{conv} \big( \Supp(g)+ \{ (\alpha,\beta, 0) \in \Z^{2n+1} |
\forall (u,v)\in \W, (u,v) \cdot (\alpha, \beta) \le 0 \}  \big)$. See
\cite[3.3]{btaka}.

\begin{clai}[Recall 2] \label{recall2}
Let $w \in \W$ and $\prec$ be an admissible order. Let $G$ be the
$\prec_w^h$-reduced standard basis of $h(I)$. Let $Q=\sum_{g \in G}
\mathrm{New}(g)$ be the Minkowski sum of the Newton polyhedra, then
\[ C_w(h(I))= \mathrm{N}_Q (\mathrm{face}_w ( Q )).\]
Here $\mathrm{N}_Q (\mathrm{face}_w (Q) )$ denotes
the normal cone in $Q$ of the face of $Q$ with respect to $w$.
\end{clai}

See Prop. 3.3.3 of \cite{btaka} for the proof and 3.3 of loc. cit. for
more details.

\begin{lem}\label{lem:prepare}
Suppose we are given $\QQ \subset \CC$ prime, $h \in \CC
\smallsetminus \QQ$ and $g \in \FDn(\CC[h^{-1}])\z$. Then 
there exists $H \in \CC \smallsetminus \QQ$ such that
for any $\PP \in V(\QQ) \smallsetminus V(H)$,
$\New((g)_\PP)=\New((g)_\QQ)$.
\end{lem}
\begin{proof}
Denote by $\W^\star$ the polar dual cone of $\W$,
$\W^\star=\{\sum_{i=1}^n \lambda_i e_i +\lambda'_i e'_i | \lambda_i,
\lambda'_i \ge 0\}$, where $e_i \in \Z^{2n}$ is the vector having $1$
in its $i$th component and $0$ for the others, and $e'_i \in \Z^{2n}$
has $-1$ at positions $i$ and $n+i$ and $0$ at the others.
Then there exists a finite subset $E(g) \subset \Supp(g)$ such that
$\New((g)_\QQ)=\mathrm{conv}(E(g)) +\W^\star$ (the
argument uses Dickson lemma, see \cite[3.3]{btaka}).

Let us write $g=\sum_\gamma \frac{c_\gamma}{h^{l_\gamma}} (x, \ddx,
z)^\gamma$, with $\gamma \in \N^{2n+1}$, $c_\gamma \in \CC$, $l_\gamma
\in \N$, and set $H=h \prod_{\gamma \in E(g)} c_\gamma(g)$.

Take $\PP \in V(\QQ) \smallsetminus V(H)$ and let us prove
$\mathrm{New}((g)_\PP)= \mathrm{New}((g)_\QQ)$.

Since $\PP \supset \QQ$, we have $\Supp((g)_\PP) \subset
\Supp((g)_\QQ)$ thus $\mathrm{New}((g)_\PP) \subset
\mathrm{New}((g)_\QQ)$. Let us prove the inverse inclusion.

We have $\Supp((g)_\QQ) \subset \mathrm{conv}(E(g)) + \W^\star$ but 
$E(g) \subset \Supp((g)_\PP)$ (because $H \notin \PP$), therefore
$\Supp((g)_\QQ) \subset \mathrm{New}((g)_\PP)$. This implies the
desired inclusion and completes the proof.
\end{proof}

\noindent
We are now ready to prove Theorem \ref{theo:main}.
\begin{proof}
The ideal $I$ is in $\C\{x,y\}[\ddx]$. We can consider a
representative of $I$ and we regard it as an ideal of
$\FDn(\CC)$ where $\CC=\O_{\C^m}(U)$ and $U$ is a polydisc in $\C^m$.
We shall prove that there exists $h \in \CC
\smallsetminus \QQ$ such that for any $\PP \in V(\QQ) \smallsetminus
V(h)$ the Gr\"obner fan $\E(h( (I)_\PP ))$ does not depend on
$\PP$. This will prove our theorem because the formal Gr\"obner fan
and the analytic Gr\"obner fan coincide.

From now on, we assume $I$ in $\FDn(\CC)$. First let us prove the
following:
\begin{clai}\label{claim:homog}
There exists $h' \in \CC \smallsetminus \QQ$ such that for a $\PP$ in
$V(\QQ) \smallsetminus V(h')$, $h((I)_\PP)= (h(I))_\PP$.
\end{clai}

Consider the order $\prec_t$ defined on $\N^{2n}$ as follows:
Take any admissible order $\prec$. Take
$t=(0,\ldots, 0, 1, \ldots, 1)$ (the number of $0$'s and $1$'s is $n$)
and refine $t$ by $\prec$ to obtain $\prec_t$. Here $t$ stands for
``total degree''. For such an order we also have a formal division 
in $\FDnk$ and the notion of standard basis (see
e.g. Castro-Jim\'enez, Granger \cite{cg}).
Therefore, we have a notion of generic standard basis
(see Bahloul \cite{jmsj}).

Now take $(\GG',h')$ a generic standard basis of $I$ for $\prec_t$.
For any $\PP \in V(\QQ) \smallsetminus V(h')$, $(\GG')_\PP$ is a
$\prec_t$-standard basis of $(I)_\PP$. Therefore, $h((I)_\PP)$ is
generated by $h( (\GG')_\PP)$ (this is well known, see
e.g. \cite[4.2]{btaka}).

Take $P$ in $I$. Consider the division modulo $\QQ$ of $P$ by $\GG'$
with $\prec_t$: $P=\sum_g q_g g + R +T$ as in
Prop. \ref{prop:divmodQ}.
By specializing to $\QQ$, $(P)_\QQ=\sum_g (q_g)_\QQ (g)_\QQ +
(R)_\QQ$, we obtain the division of $(P)_\QQ$ by $(\GG')_\QQ$ as in
Theorem \ref{theo:div}, thus $(R)_\QQ=0$ because $(\GG')_\QQ$ is a
standard basis, i.e. $R$ equals zero modulo $\QQ$, we may assume it is
zero. 
By definition of $\prec_t$, $\deg(P) \ge \deg(q_g g)$ and $\deg(P) \ge
\deg(T)$, therefore $h(P)=\sum z^{l_g} h(q_g) h(g) +z^{l} h(T)$ for
some integers $l$ and $l_g$. When specializing to $\PP$, we get
$(h(P))_\PP=\sum z^{l_g} (h(q_g))_\PP (h(g))_\PP$. We also notice that
$(h(g))_\PP = h((g)_\PP)$ (because the leading exponent which has
maximum degree is preserved after specialization to $\PP$) and that
$(h(I))_\PP$ is generated by the $(h(P))_\PP$, $P\in I$. Thus, we can
conclude that $h((I)_\PP)$ and $(h(I))_\PP$ are both generated by
$h((\GG')_\PP)$. Claim \ref{claim:homog} is proven.\\

Let us go back to the proof: in the sequel we shall use Claim
\ref{claim:homog} without any explicit mention.
Consider the (open) Gr\"obner fan of
$(I)_\QQ$: $\{ C_{w_1}(h((I)_\QQ)), \\
\ldots, C_{w_s}(h((I)_\QQ)) \}$.
For each $i$, let $(\GG_i, h_i)$ be a red.gen.s.b. of $h(I)$ for
$\prec_{w_i}^h$.

Define $h''$ as the product of all the $H$'s obtained when we apply
Lemma \ref{lem:prepare} to all the elements of the $\GG_i$'s (this
product is finite).
Now take $w \in \W$ arbitrary, put $h=h' h'' h_1 \cdots h_s$ ($h'$ comes
from Claim \ref{claim:homog})
and let us prove that for any $\PP \in V(\QQ) \smallsetminus V(h)$,
$C_w(h((I)_\PP))$ does not depend on $\PP$.

Since $\W = \cup_i C_{w_i}(h((I)_\QQ))$, there exists $i$ with $w \in
C_{w_i}(h((I)_\QQ))$. By definition, $(\GG_i)_\QQ$ is the reduced
standard basis of $h((I)_\QQ)$ for $\prec_{w_i}^h$. By Recall 1,
it is also the reduced standard basis for $\prec_w^h$. Thus by
definition, $(\GG_i, h_i)$ is a red.gen.s.b of $h(I)$ for
$\prec_w^h$. By Theorem \ref{theo:red.gen.sb}, $(\GG_i)_\PP$ is the
reduced standard basis of $h((I)_\PP)$ for $\prec_w^h$. By using
Recall 2 and lemma \ref{lem:prepare}, we can conclude that:
$C_w(h((I)_\PP))= C_w(h((I)_\QQ))$.
\end{proof}

\begin{rem*}
The progress of the proof shows that Th. \ref{theo:main} is
true for other situations (e.g. if the ideal $I$ is in $\Dn[y]$).
\end{rem*}

\end{document}